\documentclass[a4paper,11pt]{article}
\usepackage{graphicx}
\usepackage{amsfonts}
\usepackage{amsmath}
\usepackage{amssymb}
\usepackage{indentfirst}
\usepackage{titlesec}

\def\barr{\begin{array}}
\def\earr{\end{array}}
\def\bali{\begin{aligned}}
\def\eali{\end{aligned}}
\def\bearr{\begin{eqnarray}}
\def\eearr{\end{eqnarray}}

\providecommand{\play}{\displaystyle}

\providecommand{\li}{\limits}

\providecommand{\pt}{\partial}

\providecommand{\ra}{\rightarrow}

\providecommand{\da}{\downarrow}

\providecommand{\Prob}{\mathbf P}

\providecommand{\E}{\mathbf E}

\providecommand{\al}{\alpha}

\providecommand{\gm}{\gamma}

\providecommand{\Gm}{\Gamma}

\providecommand{\dt}{\delta}

\providecommand{\Dt}{\Delta}

\providecommand{\ve}{\varepsilon}

\providecommand{\tht}{\theta}

\providecommand{\lb}{\lambda}

\providecommand{\om}{\omega}

\providecommand{\R}{\mathbb R}

\providecommand{\cE}{\mathcal E}

\providecommand{\cF}{\mathcal F}

\providecommand{\cL}{\mathcal L}

\providecommand{\cN}{\mathcal N}

\providecommand{\cU}{\mathcal U}

\providecommand{\1}{\mathbf 1}

\providecommand{\contfunc}{\mathbf{C}}

\providecommand{\grad}{\nabla}

\providecommand{\vphi}{\varphi}

\providecommand{\bbar}{\bar{b}}

\begin{document}

\title{Random perturbations of dynamical systems with reflecting boundary and corresponding PDE with a small parameter}
\author{Wenqing Hu\thanks{Department of Mathematics, University of Maryland at College Park,
huwenqing@math.umd.edu} , Lucas Tcheuko\thanks{Department of
Mathematics, University of Maryland at College Park,
lucast@math.umd.edu}}
\date{}
\maketitle

\begin{abstract}
We study the asymptotic behavior of a diffusion process with small
diffusion in a domain $D$. This process is reflected at $\pt D$ with
respect to a co-normal direction pointing inside $D$. Our asymptotic
result is used to study the long time behavior of the solution of
the corresponding parabolic PDE with Neumann boundary condition.
\end{abstract}

\textit{Keywords:} PDE with a small parameter, large deviations,
Freidlin-Wentzell theory, diffusion process with reflection.

\textit{2010 Mathematics Subject Classification Numbers:} 60J60,
60F10, 60H30.

\section{Introduction}

Consider the following parabolic initial-boundary value problem

$$
\left\{\begin{array}{lr}
\displaystyle{\frac{\partial{u^\varepsilon}}{\partial{t}}=
\cL^\varepsilon u^\varepsilon \equiv
\frac{\varepsilon^2}{2}\sum\limits_{i,j=1}^{d}a_{ij}(x)\frac{\partial^2
u^\varepsilon }{\partial x_{i} \partial
x_{j}}+\sum\limits_{i=1}^{d}b^{i}(x)\frac{\partial
u^\varepsilon}{\partial
x_{i}}} \ ,  & \varepsilon > 0  \ ; \\
u^\varepsilon(x,0)=g(x) \ , &  x \in D\cup \pt D \ ; \\
\displaystyle{\frac{\partial u^\varepsilon}{\partial \gamma} (x, t)
=0} \ ,  & x \in
\partial D \ , \ t \geq 0 \ .
\end{array}
\right. \eqno(1.1)
$$

\

Here $D$ is a $d$-dimensional bounded domain in $\R^d$ with smooth
boundary $\partial D$. The initial condition $g(\bullet)$ is smooth
in $D \cup \pt D$. The matrix $a(x)=(a_{ij}(x))_{1\leq i,j\leq d}$
is positive definite. The functions $a_{ij}(x)$ are smooth and
uniformly bounded, with uniformly bounded derivatives. There is a
constant $\tht>0$ such that for any $\xi = (\xi^{1},...,\xi^{d})$ we
have $\theta^{2} |\xi|_{\R^d}^2 \leq
\sum\limits_{i,j=1}^{d}a_{ij}(x)\xi^{i}\xi^{j}\leq \theta^{-2}
|\xi|_{\R^d}^2$. The vector field $b(x)=(b^1(x),...,b^d(x))$ have
terms which are uniformly bounded, smooth in $[D]$ (here and below
$[D]$ is the closure of $D$ in Euclidean metric), and have uniformly
bounded derivatives. The vector field
$\gamma(x)=(\gamma^1(x),...,\gamma^d(x))$ is the inward
\textit{co-normal} unit vector field on $\partial D$ with respect to
the matrix $a^{-1}(x)\equiv (a^{ij}(x))_{1\leq i, j \leq
d}=(a_{ij}(x))^{-1}_{1\leq i, j \leq d}$. That is to say, for any
vector $v(x)=(v^1(x),...,v^d(x))$ tangent to $\partial D$ we have
$(\gamma, v)_{a^{-1}(x)}\equiv \sum\limits_{i,j=1}^{d}
a^{ij}(x)\gamma^{i}(x) v^{j}(x) = 0$. (Here and below
$(\gm,v)_{a^{-1}(x)}$ is the inner product with respect to the
matrix $a^{-1}(x)$. For a detailed discussion of the co-normal
condition we refer to \cite[Section 2.5]{[F red book]}.) We also
have $|\gm|_{\R^d}=1$.

Let us assume that the vector field $b(x)$ is pointing outward to
$D$ on a connected subset $\pt_1 D$ of $\pt D$, and it is pointing
inward on $\pt_2 D \equiv \pt D\backslash \pt_1 D$. (It is never
tangent to $\pt D$.) Let $\bbar(x)$ be the field coinciding with
$b(x)$ everywhere except at those points of $\pt_1 D$. At these
points $\bbar(x)$ is defined as the projection of $b(x)$ onto the
direction of the boundary. Suppose the dynamical system
$\dot{x}_t=\bbar(x_t)$ has all its $\omega$-limit sets on $\pt_1 D$.
These $\omega$-limit sets are points $O_1,...,O_l$ ($l\geq 1$).

\

Our goal in this paper is to describe the long-time behavior of the
solution $u^\varepsilon(x,t)$ of (1.1) as $\varepsilon \rightarrow
0$ and $t \rightarrow \infty$. One can relate problem (1.1) with a
certain diffusion process $X_t^\ve$ with small diffusion and
reflection with respect to $\gamma$ on $\partial D$. This process
can be described as a solution of the following stochastic
differential equation:

$$
dX_t^\varepsilon = b(X_t^\varepsilon)dt +
\varepsilon\sigma(X_t^\varepsilon)dW_t + \1_{\partial
D}(X_t^\varepsilon)\gamma(X_t^\varepsilon)d\xi_t^\varepsilon \ , \
X_0^\varepsilon =x \ , \ \xi_0^\varepsilon = 0 \ . \eqno(1.2)
$$

Here $\sigma(x)$ is a $d\times d$ matrix with smooth terms (and
bounded derivatives) that satisfies
$\sigma(x)\sigma^T(x)=\sigma^T(x)\sigma(x)=a(x)$. The function
$\1_{\partial D}(\bullet)$ is the indicator function of $\partial
D$. The processes $X_t^\varepsilon$ and $\xi_t^\varepsilon$ are
continuous time stochastic processes, adapted to the filtration
$(\cF_t)_{t\geq 0}$. They satisfy the following assumptions with
probability $1$:

(1) The process $X_t^\varepsilon \in [D]$ ;

(2) The process $\xi_t^\varepsilon$ is non-decreasing in $t$ and
increases only at $\Delta =\{t \ ; \ X_t^\varepsilon \in
\partial D\}$;

(3) The set $\Delta$ has Lebesgue measure zero.

Under these assumptions, it was proved in \cite{[Anderson-Orey]}
(also see \cite{[Watanabe]}) that such a pair of processes
$(X_t^\varepsilon, \xi_t^\varepsilon)$ exist and is unique (in the
sense of probability 1). The process $\xi_t^\varepsilon$ is called
\textit{the local time} of the process $X_t^\varepsilon$ on
$\partial D$. (We remark here that this notion of the local time for
the multidimensional diffusion process extends the classical
1-dimensional local time in \cite{[Ito-McKean]}. See
\cite{[Watanabe]} for a discussion based on SDE approach. For other
discussions of the local time for multidimensional diffusion process
we also refer to \cite{[Sato-Tanaka]} and \cite{[Sato-Ueno]}.) The
process $X_t^\varepsilon$ is a strong Markov process in $[D]$ and it
satisfies the Doeblin condition, which leads to the existence and
uniqueness of an invariant measure in $[D]$.

\

It turns out that the solution $u^\varepsilon(x,t)$ of (1.1) can be
represented as $u^\varepsilon(x,t)=\E_x g(X_t^\varepsilon)$ (see
Section 4 for details). Thus the asymptotic behavior of solution
$u^\varepsilon(x,t)$ as $\varepsilon \rightarrow 0, t\rightarrow
\infty$ is determined by the asymptotic behavior of the process
$X_t^\varepsilon$. However, the latter can be calculated using the
Freidlin-Wentzell large deviation theory (see \cite{[FW book]},
\cite{[FW 1969]}).

In Section 2 of the present paper we will give an expression of the
action functional $S^+_{0T}$ of the process $X_t^\varepsilon$. By
using the large deviation principle for the family of processes
$\{X_t^\ve\}_{\ve>0}$ we will give a description of the asymptotic
behavior of $X_t^\varepsilon$ in Section 3. Since the proof is based
on the method of \cite[Ch.6]{[FW book]} and \cite{[FW 1969]}, we
will only prove some key technical lemmas and sketch the result. In
particular, we give the algorithm on the calculation of metastable
states. Section 4 provides the corresponding result for problem
(1.1). We point out that a related question for elliptic boundary
value problems was already considered in \cite{[FZ]} (also see
\cite[Section 10.3]{[FW book]}). An example is given in Section 5.

\section{Calculation of the action functional}

In this section we give an expression of the action functional
corresponding to the large deviation principle of the process
$X_t^\varepsilon$. The main proofs and justifications of our results
are contained in \cite{[Anderson-Orey]} (also see \cite[Section
10.3]{[FW book]}), so we just summarize the results we need.

In \cite{[Anderson-Orey]}, the authors have constructed the process
$(X_t^\varepsilon, \xi_t^\varepsilon)$ corresponding to (1.2) by
first realize it in the space $\mathbb{R}^d_+$ using the following
stochastic differential equation:

$$
dY_t^\varepsilon=b(\Gamma(Y_t^\varepsilon))dt+\varepsilon
\sigma(\Gamma(Y_t^\varepsilon))dW_t \ , \ Y_0^\varepsilon=x\in
\R^d_+ \ . \eqno(2.1)
$$

Here $\Gamma : \contfunc_{[0,\infty)}(\mathbb{R}^d) \rightarrow
\contfunc_{[0,\infty)}(\mathbb{R}^d_+)$ is a functional defined by
$$\Gm(\psi_t)\equiv(\Gamma(\psi))_t\equiv\Gm_t(\psi)= (\psi_t^1-0 \wedge \inf\limits_{0\leq s\leq
t}\psi_s^1, \psi_t^2,...,\psi_t^d) \eqno(2.2)$$ for
$\psi_t=(\psi_t^1,...,\psi_t^d)\in
\contfunc_{[0,\infty)}(\mathbb{R}^d)$. It was proved in
\cite{[Anderson-Orey]} that in the case of a half space $\R^d_+$ one
can take $(X_t^\varepsilon,
\xi_t^\varepsilon)=(\Gamma(Y_t^\varepsilon),(\Gamma(Y_t^\varepsilon)-Y_t^\varepsilon)^1)$.

\

In the general case when $D$ is a bounded region in $\R^d$ with
smooth boundary one can take a finite covering of $D$ by a set of
open neighborhoods $\{\cU_1,...,\cU_N\}$. Within each $\cU_i$
($i=1,...,N$) the process can be constructed via a homeomorphism
between $\cU_i$ and $\R^d$, or between $\cU_i\cap D$ and $\R^d_+$
(when $\cU_i\cap \pt D\neq \emptyset$). In the latter case we use
the construction of the process in half space as above. By
appropriately "glue" these pieces of the trajectories together one
can construct the processes $(X_t^\ve, \xi_t^\ve)$. The process
$X_t^\ve$ is the diffusion process with reflection in $D$ and the
process $\xi_t^\ve$ is the local time on $\pt D$. For details of
this construction we refer to \cite{[Anderson-Orey]}, \cite[Section
1.6]{[F red book]}.

\

It was shown in \cite[Section 1.2]{[Anderson-Orey]} that the
corresponding action functional for the family of processes
$\{X_t^\ve\}_{\ve>0}$ as $\ve \da 0$ is given by the formula

$$ S_{0T}^+(\varphi)=\left\{
\begin{array}{l}
\play{\frac{1}{2}\int_0^T\|\dot{\varphi_s}-b(\varphi_s)-\1_{\partial
D}(\varphi_s)\omega(s)\gamma(\varphi_s)\|^2_{a^{-1}(\varphi_s)}ds} \
,
\\
 \ \ \ \ \ \ \ \ \ \ \ \ \ \ \ \text{ for } \varphi\in \contfunc_{[0,T]}([D]) \text{ absolutely
continuous }, \ \varphi_0=x \ ;
\\
+\infty \ , \ \ \ \ \ \ \ \ \ \text{ for the rest of } \varphi\in
 \contfunc_{[0,T]}([D]) \ .
\end{array}
\right. \eqno(2.3)$$

Here $$\omega(s)={\frac{(\dot{\varphi_s}-b(\varphi_s),
\gamma(\varphi_s))_{a^{-1}(\varphi_s)}}{\|\gamma(\varphi_s)\|^2_{a^{-1}(\varphi_s)}}\vee
0} \ ,$$ and $\|v\|_{a^{-1}(x)}=(v,v)^{1/2}_{a^{-1}(x)}$ for vector
$v\in \R^d$.

We have

$$\begin{array}{l} \1_{\partial
D}(\varphi_s)\omega(s)\gamma(\varphi_s)
\\=\displaystyle{ \1_{\partial D} (\varphi_s) \gamma(\varphi_s)\left(
\frac{(\dot{\varphi_s},
\gamma(\varphi_s))_{a^{-1}(\varphi_s)}}{\|\gamma(\varphi_s)\|^2_{a^{-1}(\varphi_s)}}-\frac{(b(\varphi_s),
\gamma(\varphi_s))_{a^{-1}(\varphi_s)}}{\|\gamma(\varphi_s)\|^2_{a^{-1}(\varphi_s)}}\right)
\vee 0}\\= \displaystyle{- \1_{\partial D} (\varphi_s)
\frac{\gamma(\varphi_s)}{\|\gamma(\varphi_s)\|^2_{a^{-1}(\varphi_s)}}
[0\wedge (b(\varphi_s), \gamma(\varphi_s))_{a^{-1}(\varphi_s)}]} \
\text{ for a.s. } s\in [0,T] \ .
\end{array}$$

Define

$$\bar{b}(x)=b(x)-\1_{\partial D} (x)
\frac{\gamma(x)}{\|\gamma(x)\|^2_{a^{-1}(x)}} [0\wedge (b(x),
\gamma(x))_{a^{-1}(x)}] \ . \eqno(2.4)$$

We see that $\bbar(x)$ is the field coinciding with $b(x)$
everywhere except at those points of $\pt_1 D$. (Recall that $\pt_1
D$ is the part of the boundary $\pt D$ on which $b(x)$ is pointing
outward.) At these points $\bbar(x)$ is defined as the projection of
$b(x)$ onto the direction of the boundary. The action functional for
the family of processes $\{X_t^\ve\}_{\ve>0}$ can now be formulated
as

$$S_{0T}^+(\varphi)=\left\{
\begin{array}{ll}
\displaystyle{\frac{1}{2}\int_0^T\|\dot{\varphi_s}-\bar{b}(\varphi_s)\|^2_{a^{-1}(\varphi_s)}ds
} \ ,& \text{for } \varphi\in \contfunc_{[0,T]}([D]) \text{
absolutely continuous } , \varphi_0=x \ ;\\+\infty \ , & \text{ for
the rest of } \varphi\in C_{[0,T]}([D]) \ .
\end{array}
\right. \eqno(2.5)$$

The deterministic trajectory $X_t^0$ at which the above action
functional is $0$ is also calculated in \cite{[Anderson-Orey]}. It
is given by the system $\dot{x}_t=\bar{b}(x_t)$, $x_0=x$, i.e., it
coincides with the deterministic trajectory given by the vector
field $b(x)$ everywhere except at those points of $\pt_1 D$, and at
points of $\pt_1 D$ it follows the projection of $b(x)$ onto the
direction of the boundary.

\

We formulate below the large deviation principle for the family of
processes $\{X_t^\varepsilon\}_{\ve>0}$.

\

\textbf{Theorem 2.1.} (Large deviation principle) \textit{ For the
process $X_t^\varepsilon$, we have}

\textit{(i) The set $\Phi(s)=\{\varphi \in \contfunc_{[0,T]}([D]):
S_{0T}^+(\varphi) \leq s\}$ is compact for every $s \geq 0$; }

\textit{(ii) Given $\varphi \in \contfunc_{[0,T]}([D])$. For any
$\delta
> 0$ and any $\gamma
> 0$ there exist an $\varepsilon_0 > 0$ such that for any $0 <
\varepsilon < \varepsilon_0$ we have
$$\Prob\{\rho_{0T}(X^\varepsilon, \varphi)<\delta\}\geq
\exp[-\varepsilon^{-2}(S_{0T}^+(\varphi)+\gamma)] \ , \eqno(2.6)$$
where $T
> 0$ and $\rho_{0T}(\bullet,\bullet)$ denotes the uniform distance
between functions in $\contfunc_{[0,T]}([D])$;}

\textit{(iii) For any $\delta,\gamma>0$ and any $s>0$ there exists
an $\varepsilon_0>0$ such that for any $0<\varepsilon<\varepsilon_0$
we have $$\Prob\{\rho_{0T}(X^\varepsilon, \Phi(s))\geq \delta\}\leq
\exp [ -\varepsilon^{-2}(s-\gamma)] \ , \eqno(2.7)$$ where
$\rho_{0T}(\vphi, \Phi(s))=\inf\li_{\psi\in \Phi(s)}\rho_{0T}(\vphi,
\psi)$.}

\section{Asymptotic behavior of $X_t^\varepsilon$}

\subsection{Estimates on the time to converge to $\om$-limit sets on the boundary}

We now begin our study of the asymptotic behavior of the process
$X_t^\varepsilon$. First, since the dynamical system
$\dot{x_t}=\bbar(x_t)$ does not have any $\omega$ - limit set within
$D$, we shall expect that as $\varepsilon$ is small, the
trajectories of $X_t^\varepsilon$ come to the boundary $\partial_1
D$ within finite time. (Notice that at points of $\pt_1 D$ the
vector field $b(x)$ is pointing outward and at points of $\pt_2 D$
it is pointing inward. Therefore the deterministic trajectory
$X_t^0$ will not come to $\pt_2 D$.)

For any $x,y \in  [D]$, we define
$$V^+(x,y)=\inf\limits_{\varphi \in \contfunc_{[0,T]}
([D])}\{S_{0T}^+(\varphi), \varphi_0=x, \varphi_T=y, \varphi_t \in D
\cup \partial D, 0 \leq t \leq T<\infty\} \ .$$

Recall that the dynamical system $\dot{x}_t=\bbar(x_t)$ has all its
$\omega$-limit sets on $\pt_1 D$. These $\omega$-limit sets are
points $O_1,...,O_l$ ($l\geq 1$). Let us suppose, that for any $x$
and $y$ in $[D]$, $x\neq y$ we have at least one of $V^+(x,y)$ and
$V^+(y,x)$ being $>0$.

For each $O_i, i=1,2,...,l$, by an $\alpha$-\textit{neighborhood}
$\mathcal{E}_\alpha(O_i)$ of $O_i$, we refer to the intersection of
$D$ with an open ball having center $O_i$ and radius $\alpha>0$. We
use the symbol $\partial\mathcal{E}_\alpha(O_i)$ to mean the
intersection of $[D]$ with the boundary of the open $\alpha$-ball
centered at $O_i$. We call $\pt\cE_\al(O_i)$ the \textit{boundary}
of the $\alpha$-neighborhood of $O_i$. Let us choose $\alpha> 0$
such that the $\alpha$-neighborhoods $\mathcal{E}_{\alpha}(O_i)$ for
all $O_i, i=1,2,...,l$, does not intersect each other. We now prove
the following:

\

\textbf{Theorem 3.1.} \textit{There exist positive constants $c$ and
$T_0$ such that for all sufficiently small $\varepsilon > 0$ and any
$x \in [D], X_0^{\varepsilon}=x$ we have
$$\Prob_x\{\zeta_{\alpha} > T\} \leq \exp[-\varepsilon^{-2}c(T-T_0)] \ ,$$
where $\zeta_\alpha = \inf\{t: X_t^{\varepsilon}\in
\bigcup\limits_{i=1}^{l} [\mathcal{E}_{\alpha}(O_i)]\}$ \ .}

\

\textbf{Proof.}  We consider the dynamical system
$\dot{x_t}=\bar{b}(x_t)$ on the whole domain $D\cup\partial D$,
where vector field $\bar{b}(x)$ is defined as before. Since system
$\dot{x_t}=\bbar(x_t)$ does not have any $\omega$ - limit set in $D$
and $(b(x),\gamma(x))_{a^{-1}(x)}|_{\partial_2 D}
>0 $, we can say that the time $T_1(x)$ that the trajectory $x_t(x)$
spends until reaching $\partial_1 D$ is finite (if $x\in \partial_1
D$, let $T_1(x)=0$). Let $y(x)$ be the point where trajectory first
hits $\partial_1 D$. Starting from $y(x)$, the time
$T_2(y(x),\alpha)=T_2(x,\alpha)$ that the trajectory of system
$\dot{x_t}=\bar{b}(x_t)$ on $\partial_1 D$ spend to come into
$\bigcup\limits_{i=1}^{l} [\mathcal{E}_{\frac{\alpha}{2}}(O_i)]$ is
also finite (as is the same, if $y(x) \in \bigcup\limits_{i=1}^l
[\mathcal{E}_{\frac{\alpha}{2}}(O_i)]$, then $T_2(x,\alpha)=0$). The
function $T(x,\alpha)=T_1(x)+T_2(x,\alpha)$ is upper semi-continuous
in $x$ (i.e., for $x, x_0\in [D]$ we have $\overline{\lim}_{x\ra
x_0}T(x,\al)\leq T(x_0,\al)$) because $x_t(x)$ depends continuously
on $x$. Thus there exists $T_0=\max\limits_{x \in [D]}T(x,\alpha) <
\infty$. The set of functions in $\contfunc_{[0,T_0]}([D])$ assuming
their values in $[D] \backslash \left(\bigcup\limits_{i=1}^l
\mathcal{E}_{\frac{\al}{2}}(O_i)\right)$ is closed and thus
$S_{0T_0}^+$ attains a minimum $A$ on this set. Taking into account
the construction of $T_0$ and the form of $S_{0T_0}^+$ in (2.5), we
see that $A>0$. Let $0 < \delta < \frac{\alpha}{2}$. Let
$\Phi_x(A/2)=\{\varphi \in \contfunc_{[0,T_0]}(D), \varphi_0=x,
S_{0T_0}^+(\vphi) \leq A/2\}$. We see that trajectories for which
$\zeta_\alpha > T_0$ are at a distance $\geq \delta$ from
$\Phi_x(A/2)$. Thus by the part (iii) of the large deviation
principle we have
$$\Prob_x\{\zeta_\alpha > T_0\} \leq \exp[-\varepsilon^{-2}(A/2-\gamma)]$$
for some $0<\gamma<A/2$.

Thus by strong Markov property,

$$
\begin{array}{ll}
\Prob_x\{\zeta_\alpha > (n+1)T_0\} & =\E_x[\zeta_\alpha > nT_0;
\Prob_{X_{nT_0}^\varepsilon}\{\zeta_\alpha > T_0\}]\\
& \leq \Prob_x\{\zeta_\alpha > nT_0\}
\exp[-\varepsilon^{-2}(A/2-\gamma)] \ .
\end{array}
$$
So by induction we see that

$$
\begin{array}{ll}
\Prob_x\{\zeta_\alpha > T\} & \displaystyle{\leq
\Prob_x\{\zeta_\alpha
> \left[\frac{T}{T_0}\right]T_0\}} \\ & \displaystyle{\leq
\exp\{-\varepsilon^{-2}(A/2-\gamma)\left[\frac{T}{T_0}\right]\}} \\
& \displaystyle{\leq
\exp\{-\varepsilon^{-2}\left(\frac{T}{T_0}-1\right)(A/2-\gamma)\}} \
.
\end{array}
$$
Putting $c=\displaystyle{\frac{A/2-\gamma}{T_0}}$, we get as
desired. $\square$

\

\subsection{Transition probabilities between neighborhoods of the $O_i$'s}

In this section we study the asymptotic transition probabilities
between neighborhoods of the $\om$-limit sets $\{O_1,...,O_l\}$. We
first provide several auxiliary lemmas.

\

\textbf{Lemma 3.1.} \textit{There exists a constant $L>0$ such that
for any $x, y\in [D]$ sufficiently close to each other, there exists
a function $\varphi \in \contfunc_{[0,T]}([D]), \varphi_0=x,
\varphi_T=y$, such that we have
$S_{0T}^+(\varphi)<L\cdot|x-y|_{\R^d}$.}

\

\textbf{Proof.} Let $x$ and $y$ be so close to each other that they
can be covered by one coordinate chart $U$. Let this coordinate
chart correspond to a coordinate function $u: U \ra \R^d (\text{ or
} \R^d_+)$. The function $u$ is smooth with bounded derivatives. Let
us take $T=|x-y|_{\R^d}$,
$$\varphi_t= u^{-1}\left( u(x)+\frac{t}{T}(u(y)-u(x)) \right) \ .$$

We have, for some constant $M>0$,

$$
\begin{array}{ll}
S^+_{0T}(\varphi_s)
&\displaystyle{=\frac{1}{2}\int_0^T\|\dot{\varphi_s}-\bar{b}(\varphi_s)\|^2_{a^{-1}(\varphi_s)}ds}
\\ &\displaystyle{=\frac{1}{2}\int_0^T\sum\limits_{i,j=1}^d
a^{-1}_{ij}(\varphi_s)(\dot{\varphi}_s^i-\bar{b}^i(\varphi_s))(\dot{\varphi}_s^j-\bar{b}^j(\varphi_s))ds}\\&\displaystyle{\leq
\frac{\theta^2}{2}\int_0^T|\dot{\varphi_s}-\bar{b}(\varphi_s)|_{\R^d}^2ds}
\\&=\displaystyle{\frac{\theta^2}{2}\int_0^T\left|\frac{1}{T}(u^{-1})'\left(u(x)+\frac{t}{T}(u(y)-u(x))\right)
\circ(u(y)-u(x))-\bar{b}(\varphi_s)\right|_{\R^d}^2ds}\\
&\displaystyle{\leq \theta^2M\left\{\frac{1}{T}|y-x|_{\R^d}^2 +T
\max\limits_{x\in [D]}|\bar{b}(x)|_{\R^d}^2\right\}} \ .
\end{array}
$$

Taking into account that $T=|y-x|_{\R^d}$, we are done. $\square$

\

\textbf{Lemma 3.2.} \textit{For any $\gamma > 0$ and any compact
subset $K\subseteq[D]$ there exists $T_0$ such that for any $x,y \in
K$ there exists a function $\varphi_t, 0 \leq t \leq T,
\varphi_0=x,\varphi_T=y, T \leq T_0$ such that
$S_{0T}^+(\varphi)\leq V^+(x,y)+\gamma$.}

\

\textbf{Proof.} We choose a finite $\delta$-net $\{x_i\}$ of points
in $K$; we connect them with curves at which the action functional
assumes values differing from the infimum by less than
$\frac{\delta\gm}{2}$ and complete them with end sections using
Lemma 3.1: from $x$ to a point $x_i$ near $x$ and then from $x_i$ to
a point $x_j$ near $y$, and from $x_j$ to $y$. By choosing $\dt$
small enough we get as desired. $\square$

\

We define

\

$$\widetilde{V}^+(O_i,O_j)=\inf\limits_{\varphi \in
\contfunc_{[0,T]}([D])}\{S_{0T}^+(\varphi): \varphi_0=O_i,
\varphi_T=O_j, \varphi_t \in [D] \setminus \bigcup\limits_{s\neq
i,j}\{O_s\} , 0<t<T\} \ .$$

A "$\widetilde{V}^+(O_i,O_j)$ version" of the above Lemma can proved
similarly: one can take the curve $\vphi$ in such a way that it
avoids $\bigcup_{s\neq i,j}\{O_s\}$ and such that $S_{0T}(\vphi)\leq
\widetilde{V}^+(O_i,O_j)+\gm$. We omit the proof.

\

Let constant $\rho_0>0$ be small. Let constant $0 < \rho_1 <
\rho_0$. We denote by $C$ the set $D\cup
\partial D$ from which we delete the $\rho_0$-neighborhoods of the
$O_i, i=1,2,...,l$; by $\Gamma_i$ the boundaries of the
$\rho_0$-neighborhoods of $O_i$: $\Gamma_i=\partial
\mathcal{E}_{\rho_0}(O_i)$; by $g_i$ the $\rho_1$-neighborhoods of
the $O_i$, and by $g$ the union of all the $g_i$.

We introduce the following random times $\tau_0=0, \sigma_n=\inf\{t
\geq \tau_n, X_t^\varepsilon \in C\}, \tau_n=\inf\{t \geq
\sigma_{n-1}, X_t^\varepsilon \in \partial g\}$. We consider the
Markov chain $Z_n=X_{\tau_n}^\varepsilon$ for $n \geq 0$. We see
that from $n=1$ on $Z_n \in
\partial g$. Also, $X_{\sigma_0}^\varepsilon$ can be any point of
$C$, all the following $X_{\sigma_n}^\varepsilon$ belong to one of
the $\Gamma_i$'s. The chain never stops.

We are now ready to prove:

\

\textbf{Theorem 3.2.} \textit{For any $\gamma >0$ there exists
$\rho_0>0$ (which can be chosen arbitrary small) such that for any
$\rho_2$, $0<\rho_2<\rho_0$, there exists $\rho_1$,
$0<\rho_1<\rho_2$ such that for all $x$ in the $\rho_2$-neighborhood
of $O_i (i=1,...,l)$ the one-step transition probabilities of $Z_n$,
$Z_0=x$ satisfy the inequality
$$\exp[-\varepsilon^{-2}(\widetilde{V}^+(O_i,O_j)+\gamma)]
\leq \Prob(x,\partial g_j) \leq
\exp[-\varepsilon^{-2}(\widetilde{V}^+(O_i,O_j)-\gamma)]$$}
\textit{for some $0<\varepsilon<\varepsilon_0$.}

\

\textbf{Proof.} We can assume $\widetilde{V}^+(O_i,O_j) < \infty$.
Set $\widetilde{V}^+_0=
\max\limits_{i,j=1,2,...,l}\widetilde{V}^+(O_i,O_j)$. Choose
positive $\rho_0$ small enough. For every pair $O_i, O_j$ for which
$\widetilde{V}^+(O_i,O_j) < \infty$ we choose a function
$\varphi_t^{i,j}\in \contfunc_{[0,T]}([D]), 0\leq t \leq
T=T(O_i,O_j)$, such that $\varphi_0^{i,j}=O_i, \varphi_T^{i,j}=O_j,
\varphi_t^{i,j}$ does not touch $\bigcup\limits_{s \neq
i,j}\{O_s\}$, and such that (by Lemma 3.2)
$$S^+_{0T}(\varphi^{i,j})\leq \widetilde{V}^+(O_i,O_j)+0.5\gamma \ .$$

We choose positive $\rho_1$ smaller than $\frac{\rho_0}{2}$,
$\rho_2$ and $$\frac{1}{2}\min\{\rho(\varphi_t^{i,j},\bigcup_{s \neq
i,j}\{O_s\}): 0\leq t\leq T, i,j=1,2,...l\} \ .$$

For every $x$ in a $\rho_2$-neighborhood of $O_i$ we take a curve
connecting $x$ with $O_i$ and for which the value of $S^+$ does not
exceed $0.3\gamma$ (by Lemma 3.1). We combine this curve with the
curve $\varphi_t^{i,j}$ and obtain a function $\varphi_t, 0\leq
t\leq T, \varphi_0=x, \varphi_T=O_j$ (with a possible small change
of $T$ from $T=T(O_i,O_j)$) such that
$$S^+_{0T}(\varphi)\leq \widetilde{V}^+(O_i,O_j)+0.8\gamma \ .$$

From Lemma 3.2 we choose a $T_0 \geq T$, and extend the curve
$\varphi_t$ to $T\leq t\leq T_0$ by using a trajectory of the
dynamical system $\dot{x_t}=\bar{b}(x_t)$ on $\partial_1 D$, without
changing the value of $S^+_{0T_0}(\varphi)$ from that of
$S^+_{0T}(\vphi)$. We choose positive $\delta$ less than $\rho_1,
\rho_0-\rho_2$. For a trajectory of $X_t^\varepsilon$ starting from
$x$, passing at a distance from $\varphi_t$ smaller than $\delta$
for $0\leq t\leq T_0$, it must intersect with $\Gamma_i$ and reaches
the $\delta$-neighborhood of $O_j$ without getting closer than
$\delta$ from any of the other $O_s, s\neq i,j$. Moreover,
$X_{\tau_1}^\varepsilon \in
\partial g_j$, thus $$\Prob(x, \partial g_j)\geq
\Prob_x\{\rho_{0T_0}(X^\varepsilon,\varphi)<\delta\}\geq
\exp[-\varepsilon^{-2}(S^+_{0T_0}(\varphi)+0.1\gamma)]>
\exp[-\varepsilon^{-2}(\widetilde{V}^+(O_i,O_j)+\gamma)]\ .$$

Now we turn to the proof of the upper estimates. For any curve
$\varphi_t, 0\leq t \leq T$ beginning at $x$, touching the
$\delta$-neighborhood of $\partial g_j$, not touching any of the
$O_s, s \neq i,j$, we have $$S^+_{0T}(\varphi)\geq
\widetilde{V}^+(O_i,O_j)-0.7\gamma \ .$$

We use Theorem 3.1 to choose $T_1$ such that for any $x \in
[D]\setminus g$ we have $\Prob^\varepsilon_x\{\tau_1>T_1\} \leq
\exp(-\varepsilon^{-2}V_0^+)$ for some $V_0^+>0$.

Any trajectory $X_t^\varepsilon$ beginning at $x$ and being in
$\partial g_j$ at time $\tau_1$ either spends time $T_1$ without
touching $\pt g$ or reaches $\partial g_j$ over time $T_1$, in this
case
$$\rho_{0T_1}(X^\varepsilon,
\Phi_x(\widetilde{V}^+(O_i,O_j)-0.7\gamma))\geq \delta \ .$$

Therefore we have

$$\begin{array}{ll}
\Prob^\varepsilon_x\{X_{\tau_1}^\varepsilon \in
\partial g_j\} & \leq \Prob^\varepsilon_x\{\tau_1>T_1\} + \Prob^\varepsilon_x\{\rho_{0T_1}(X^\varepsilon,
\Phi_x(\widetilde{V}^+(O_i,O_j)-0.7\gamma))\geq \delta\} \\& \leq
\exp(-\varepsilon^{-2}V^+_0)+\exp[-\varepsilon^{-2}(\widetilde{V}^+(O_i,O_j)-0.9\gamma)]
\\& \leq \exp[-\varepsilon^{-2}(\widetilde{V}^+(O_i,O_j)-\gamma)]\end{array}$$
for sufficiently small $\varepsilon$. $\square$

\

In an exactly similar way one can also formulate the estimate on
transition probability based on the quantities

$$\widetilde{V}^+(x,O_j)=\inf\limits_{\varphi \in
\contfunc_{[0,T]}([D])}\{S_{0T}^+(\varphi): \varphi_0=x,
\varphi_T=O_j, \varphi_t \in [D] \setminus \bigcup\limits_{s\neq
j}\{O_s\} , 0<t<T\} \ .$$

We have

\

\textbf{Theorem 3.3.} \textit{For any $\gamma >0$ there exists
$\rho_0>0$ (which can be chosen arbitrary small) such that for any
$\rho_2$, $0<\rho_2<\rho_0$, there exists $\rho_1$,
$0<\rho_1<\rho_2$ such that for all $x$ outside the
$\rho_2$-neighborhood of $O_i (i=1,...,l)$ the one-step transition
probabilities of $Z_n$, $Z_0=x$ satisfy the inequality
$$\exp[-\varepsilon^{-2}(\widetilde{V}^+(x,O_j)+\gamma)]
\leq \Prob(x,\partial g_j) \leq
\exp[-\varepsilon^{-2}(\widetilde{V}^+(x,O_j)-\gamma)]$$}
\textit{for some $0<\varepsilon<\varepsilon_0$.}

\subsection{The invariant measure of $X_t^\ve$; sublimiting distribution}

In this section we study the invariant measure of the process
$X_t^\varepsilon$. Based on the estimates on transition
probabilities given above, the proof of the asymptotic result is the
same as that of \cite[Ch.6]{[FW book]} and \cite{[FW 1969]}. Let us
formulate and prove two more technical lemmas, after which the rest
of the proof is just a study of Markov chains on graphs. The latter
part will be omitted since it is the same as \cite[Ch.6]{[FW book]}
and \cite{[FW 1969]}.

\

\textbf{Lemma 3.3.} \textit{For $i\in\{1,2,...,l\}$, define
$$\tau_{\cE_\dt(O_i)} = \inf \{t, X_0^\varepsilon=x, X_t^\varepsilon\in \partial
\cE_\dt(O_i)\} \ .$$ For any $\gamma>0$, there exist $\delta>0$ such
that for all sufficiently small $\varepsilon$ and $x\in
\cE_\dt(O_i)$ we have
$$\E_x^\varepsilon \tau_{\cE_\dt(O_i)} < \exp(\gamma \varepsilon^{-2}) \ .$$ }

\textbf{Proof.} Choose point $z\in D$ close to $O_i$. Put
$\delta=\frac{|z-O_i|}{2}$. Connect $x$ with $O_i$ and $O_i$ with
$z$ with the values of $S^+$ not exceeding $\frac{\gamma}{4}$ and
$\frac{\gamma}{2}$, the resulting function is called
$\widetilde{\varphi}_t$. The length of the time interval of
$\widetilde{\varphi}_t$ is uniformly bounded by $T_0$ for all $x\in
G$. We extend $\widetilde{\varphi}_t$ up to $T_0$ by using a
trajectory of $\dot{x}_t=\bbar(x_t)$ in $D\cup\partial D$ without
making $S^+$ larger.

Now we have for $x\in \cE_\dt(O_i)$,
$$\Prob_x^\varepsilon\{\tau_{\cE_\dt(O_i)}<T_0\}
\geq
\Prob_x^\varepsilon\{\rho_{0T_0}(X^\varepsilon,\widetilde{\varphi})<\delta\}
\geq \exp(-0.9\gamma \varepsilon^{-2}) \ .$$

Using the Markov property we see that
$$\Prob_x^\varepsilon\{\tau_{\cE_{\dt}(O_i)}\geq nT_0\}\leq [1-\exp(-0.9\gamma \varepsilon^{-2})]^n \ .$$

This yields $$\E_x^\varepsilon \tau_{\cE_\dt(O_i)}\leq
T_0\sum\limits_{n=0}^\infty[1-\exp(-0.9\gamma
\varepsilon^{-2})]^n=T_0\exp(0.9\gamma\varepsilon^{-2}) \ .$$

Sacrificing $0.1\gamma$ in order to get rid of $T_0$ we get the
desired result. $\square$

\

\textbf{Lemma 3.4.} \textit{For any $\gamma>0$ there exist
$\rho_1>0$ such that for all sufficiently small $\varepsilon$ and $y
\in
\partial g_i$ we have
$$\E_y^\varepsilon\int_0^{\sigma_0}\chi_{g_i}(X_t^\varepsilon)dt>\exp(-\gamma \varepsilon^{-2}) \ .$$}

\textbf{Proof.} Choose $\rho_1$ small. We connect $y\in
\partial g_i$ with $O_i$ using a curve $\varphi_t$, extend it using
the trajectory of $\dot{x}_t=\bar{b}(x_t)$ on $\partial D$ till
first exit time $\sigma_0$ from $\mathcal{E}_{\rho_0}(O_i)$, with
corresponding $S^+$ less than $0.5\gamma$. All the trajectories at a
distance less than $\frac{\rho_1}{2}$ spends a time at least $t_0>0$
within $g_i$, uniformly for all $y\in \partial g_i$. The probability
of all such trajectories is no less than
$\exp(-0.9\gamma\varepsilon^{-2})$. Thus the expected value is no
less than $t_0\exp(-0.9\gamma\varepsilon^{-2})$. By sacrificing
$0.1\gamma$ we can get rid of $t_0$. $\square$

\

The rest of this section is devoted to the description of the
algorithm for the calculation of the invariant measure and the
metastable states. The proof we shall omit here follows
\cite[Ch.6]{[FW book]} and \cite{[FW 1969]}.

Let $L$ be a finite set (in our case $L=\{1,2,...,l\}$), whose
elements are denoted by letters $i,j,k,m,n$, etc. Let a subset $W$
be selected in $L$. A graph consisting of arrows $m \ra n$ ($m \in L
\backslash W, n \in L, n \neq m$) is called a $W$-graph if it
satisfies the following conditions:

(1) every point $m \in L \backslash W$ is the initial point of
exactly one arrow;

(2) there are no cycles in the graph.

Intuitively, a $W$-graph is a graph consisting of arrows starting
from each point $m \in L \backslash W$, and going along a sequence
of arrows leading to some point $n \in W$.

The set of $W$-graphs is denoted by $G(W)$. We shall use the letter
$g$ to denote graphs.

Let $W(O_i)=\min\limits_{g\in G\{i\}}\sum_{(m\rightarrow n)\in
g}\widetilde{V}^+(O_m,O_n)$. It can be proved that
$$W(O_i)=\min\limits_{g\in G\{i\}}\sum_{(m\rightarrow n)\in
g}V^+(O_m,O_n) \ .$$

We have

\

\textbf{Theorem 3.4.} \textit{Let $\mu^\varepsilon$ be the
normalized invariant measure of the process $X_t^\varepsilon$. Then
for any $\gamma>0$ there exists $\rho_1>0$ such that we have
$$\exp[-\varepsilon^{-2}(W(O_i)-\min\limits_{i}W(O_i)+\gamma)]\leq\mu^\varepsilon(g_i)
\leq \exp[-\varepsilon^{-2}(W(O_i)-\min\limits_{i}W(O_i)-\gamma)]$$
for sufficiently small $\varepsilon>0$.}

\

We shall say that a set $N\subset [D]$ is \textit{stable} if for any
$x\in N$, $y \not\in N$ we have $V^+(x,y)>0$. One can show that for
an unstable $O_j$ ($j=1,...,l$) there exist a stable $O_i$ ($i\neq j
, i=1,...,l$) such that $V^+(O_i, O_j)=0$.

\

\textbf{Theorem 3.5.} \textit{For $x\in [D]$ set}

$$W(x)=\min[W(O_i)+V(O_i,x)] \ ,$$
\textit{where the minimum can be taken over either all of
$O_1,...,O_l$ or only stable ones. Let $\mu^\varepsilon$ be the
normalized invariant measure of the process $X_t^\varepsilon$. Then
for any $\gamma>0$ there exists $\bar{\rho}>0$ such that for any
$0<\rho<\bar{\rho}$ we have
$$\exp[-\varepsilon^{-2}(W(x)-\min\limits_{i}W(O_i)+\gamma)]\leq\mu^\varepsilon(\cE_{\rho}(x))
\leq \exp[-\varepsilon^{-2}(W(x)-\min\limits_{i}W(O_i)-\gamma)]$$
for sufficiently small $\varepsilon>0$.}

\

Here $\cE_{\rho}(x)$ is a $\rho$-neighborhood of $x$.

\

The above two theorems roughly say that as first $t\ra \infty$ and
then $\ve\ra 0$, the process $X_t^\ve$ will be situated in one of
the $O_i$'s which minimizes the values of $W(O_i)$ (it can be
calculated either via all $O_1,...,O_l$ or only via the stable
ones). In generic case, when $\min\limits_i W(O_i)$ is attained at
some unique point $i$, we have for any $\delta
>0$,
$$\lim\limits_{\varepsilon \rightarrow 0}
\lim\limits_{t\rightarrow\infty}\Prob_x^\varepsilon\{|X_t^\varepsilon-O_i|>\delta\}=0
\ . \eqno(3.1)$$

\

A natural question is that how the limiting distribution behaves
when we take the limit in a coordinated way, i.e. take
$\varepsilon\rightarrow 0$ and
$t=t(\varepsilon^{-2})\rightarrow\infty$. This is the problem of
metastability and sublimiting distributions (see \cite{[F
sublimiting]}). Let us assume that
$T=T(\varepsilon)\asymp\exp(\frac{\lambda}{\varepsilon^2})$ and
 we consider
$\lim\limits_{\varepsilon\rightarrow0}\Prob_x^\varepsilon\{X_{T(\varepsilon)}^\varepsilon\in
\Gamma\}$. In the generic case one can define a function
$K^*(x,\lambda)\in \{1,2,...,l\}$ such that
$$\lim\limits_{\varepsilon\rightarrow0}
\Prob_x^\varepsilon\{|X_{T(\varepsilon)}^\varepsilon-
O_{K^*(x,\lambda)}|>\delta\}=0 \eqno(3.2)$$ for any $\delta>0$.

The algorithm to determine $K^*(x,\lambda)$ is as follows. First we
consider for each $O_i$ (the rank $0$ cycle) the "next" most
probable $\om$-limit set $\mathcal{N}(O_i)$ that we are going to
jump to. Continuing this determination of "next" states we form the
rank $1$ cycle
$O_i\ra\mathcal{N}(O_i)\ra\mathcal{N}^2(O_i)\ra...\ra\mathcal{N}^{m_{\{i\}}^1}(O_i)$.
We stop once we get a repetition $\cN(\cN^{m_{\{i\}}^1}(O_i))=O_i$.
Cycles generated by distinct initial points $i\in \{1,...,l\}$
either do not intersect each other or coincide: in the latter case
the cycle order on them is one and the same.

We continue by recurrence. Let the cycles of rank $(k-1)$ be
$\pi_1^{k-1},...,\pi_{n_{k-1}}^{k-1}$. Starting from each
$(k-1)$-cycle $\pi_{i}^{k-1}$ one can determine the "next" most
probable $(k-1)$-cycle $\pi_{\cN(\pi_i^{k-1})}^{k-1}$ that we will
first jump to. Continuing this determination we form a rank $k$
cycle $\pi_i^{k-1}\ra \pi_{\cN(\pi_i^{k-1})}^{k-1}\ra ... \ra
\pi_{\cN^{m_{\pi_i^{k-1}}}(\pi_i^{k-1})}^{k-1}$. We stop once we get
a repetition $\cN(\cN^{m_{\pi_i^{k-1}}}(\pi_i^{k-1}))=\pi_i^{k-1}$.
Cycles of rank $k$ generated by distinct cycles of rank $k-1$ either
do not intersect each other or coincide.

In this way we can continue until the last cycle which is the whole
of $\{O_1,...,O_l\}$. The metastable states are determined by the
timescale of the cycles that we traverse.

Let us be more precise. Starting from a cycle $\pi$, to determine
the "next" cycle $\cN(\pi)$ that we first jump to, we calculate

$$A(\pi)=\min\li_{g\in G(L\backslash \pi)}
\sum\li_{(m \ra n)\in g}V^+(O_m,O_n) \ . \eqno(3.3)$$

Here $L=\{1,2,...,l\}$. The minimum of the above expression
determines a $L\backslash \pi$ graph consisting of chains of arrows
leading to the first state in $L \backslash \pi$ we jump to.

We put

$$C(\pi)=A(\pi)-\min\li_{i \in \pi}\min\li_{g\in G_\pi\{i\}}
\sum\li_{(m\ra n)\in g}V^+(O_m,O_n) \ . \eqno(3.4)$$ Here
$G_{\pi}\{i\}$ is the set of $\{i\}$-graphs restricted to $\pi$.
Then the asymptotic exit time from $\pi$ is of order
$\asymp\exp\left(\frac{C(\pi)}{\ve^2}\right)$.

Starting from $i=i(x)$ (which is the label for the first equilibrium
among $O_1,...,O_l$ that we approach in finite time, starting from
$x$), let $\pi$, $\pi'$,..., $\pi^{(s)}$ be cycles of next to the
last rank, unified into the last cycle, which exhausts
$\{1,2,...,l\}$. If the constant $\lb$ is greater than
$C(\pi),C(\pi'),...,C(\pi^{(s)})$, then over time of order
$\exp(\lb\ve^{-2})$ the process can traverse all these cycles many
times (and all cycles of smaller rank inside them) and the limiting
distribution is concentrated on that one of the cycles for which
$C(\pi), C(\pi'),...,C(\pi^{(s)})$ is the greatest. Within this
cycle, it is concentrated on that one of the subcycles for which the
corresponding constant $C(\bullet)$ in (3.4) is the greatest
possible, and so on up to points (one point in the generic case)
$O_{K^*(x,\lb)}$. This point $O_{K^*(x,\lb)}$ is the metastable
state in (3.2).

\section{Application to PDE}

The solution of (1.1) can be represented through process (1.2) by
the formula $u^\varepsilon(x,t)=\E_x g(X_t^\varepsilon)$. This is an
immediate consequence of the following generalized It\^{o}'s
formula:

\

\textbf{Lemma 4.1.} \textit{Assume process $(X_t^\varepsilon,
\xi_t^\varepsilon)$ is given by (1.2), $X_0^\varepsilon=x$. Let
$u(x,t)$ be of class $\contfunc^{2,1}(\R^d\times \R_+)$ with uniform
bounded derivatives up to the second order in $x$ and up to the
first order in $t$. Then we have}

\begin{equation*}
\begin{array}{l}
u(X_t^\varepsilon,t)-u(x,0)\\=\displaystyle{\int_0^t\left(\frac{\partial}{\partial
s} + \cL^\varepsilon\right)u(X_s^\varepsilon,s)ds +\int_0^t{\nabla
u(X_s^\varepsilon,s)\cdot \gamma(X_s^\varepsilon)}d\xi_s^\varepsilon
+ \int_0^t{\nabla u(X_s^\varepsilon,s)\cdot
\sigma(X_s^\varepsilon)}dW_s} \ .
\end{array}
\end{equation*}

\

For a proof of this theorem see \cite[Section 3]{[Gikman-Skrokhod]}.

\

Our answer to the problem (1.1) is

\

\textbf{Theorem 4.1.} \textit{Under all our assumptions, in generic
case, for $T(\varepsilon)\asymp \exp
(\frac{\lambda}{\varepsilon^2})$, we have}

$$\lim\limits_{\varepsilon \rightarrow 0}u^\varepsilon(x,T(\varepsilon))= g (O_{K^*(x,\lambda)}) \
,$$\textit{where $K^*(x,\lambda)$ is defined as in Section 3.3.}

\

\section{Example}

Consider an example. Let the domain $D$ be a unit disk
$B(1)=\{(y_1,y_2); y_1^2+y_2^2<1\}$ in $\R^2$. Let the smooth vector
field $b_y(y_1,y_2)$ be given such that
$\bbar_y(y_1,y_2)=(\bbar_{y_1}(y_1,y_2),\bbar_{y_2}(y_1,y_2))$ is as
in Fig.1. We consider the problem

\begin{figure}
\centering
\includegraphics[height=7cm, width=7cm , bb=92 18 312 241]{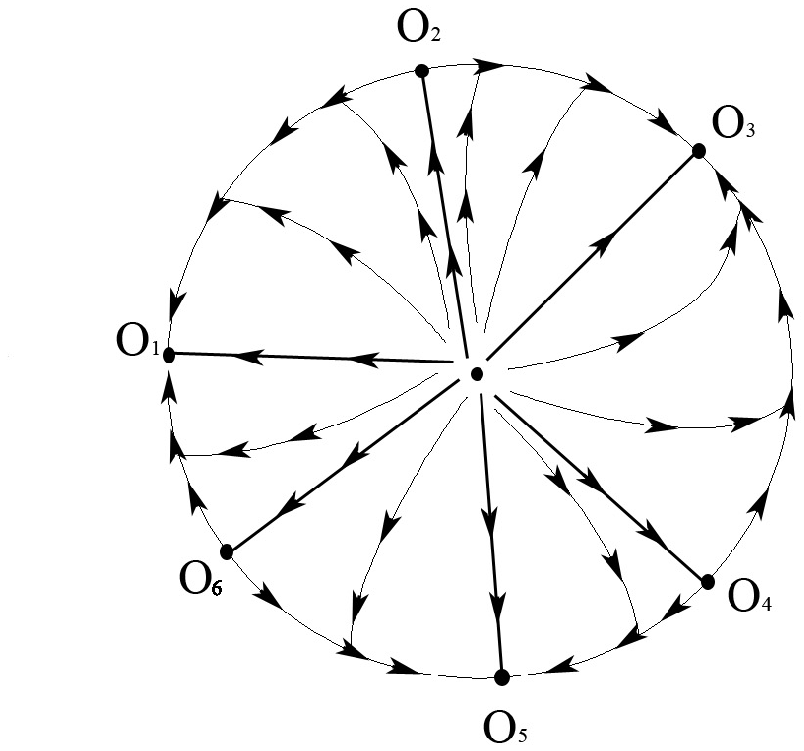}
\caption{An example}
\end{figure}

$$
\left\{\begin{array}{lr}
\displaystyle{\frac{\partial{u^\varepsilon}(y_1,y_2,t)}{\partial{t}}=
\dfrac{\ve^2}{2}\Dt_{y_1,y_2} u^\ve(y_1,y_2,t)+b_y(y_1,y_2)\cdot \grad u^\ve(y_1,y_2,t)} \ ,  & \varepsilon > 0  \ ; \\
u^\varepsilon(y_1,y_2,0)=g(y_1,y_2) \ , &  y_1^2+y_2^2\leq 1 \ ; \\
\displaystyle{\frac{\partial u^\varepsilon}{\partial r} (y_1, y_2,
t) =0} \ ,  & y_1^2+y_2^2=1 \ , \ t \geq 0 \ .
\end{array}
\right. \eqno(5.1)
$$

Here $\dfrac{\pt}{\pt r}$ is the derivative with respect to the
inward unit normal. The action functional takes the form

$$S_{0T}^+(\varphi)=\left\{
\begin{array}{ll}
\displaystyle{\frac{1}{2}\int_0^T|\dot{\varphi_s}-\bar{b}_y(\varphi_s)|^2_{\R^2}ds
} \ ,& \varphi\in \contfunc_{[0,T]}([D]) \text{ absolutely
continuous } , \varphi_0=x \ ;\\+\infty \ , & \text{ for the rest of
} \varphi\in \contfunc_{[0,T]}([D]) \ .
\end{array}
\right. \eqno(5.2)$$

We calculate the "quasi-potential" using (5.2)

$$V^+(x,y)=\inf\limits_{\varphi \in \contfunc_{[0,T]}
([D])}\{S_{0T}^+(\varphi), \varphi_0=x, \varphi_T=y, \varphi_t \in D
\cup \partial D, 0 \leq t \leq T<\infty\} \ .$$

The $\om$-limit sets of the dynamical system
$\dot{x}_t=\bbar_y(x_t)$ are the zeros of the vector field
$\bbar_y(x)$ on $\pt D=S^1$. (And also the origin but it is unstable
so that we neglect it.) In Fig.1 the points $O_1$, $O_3$ and $O_5$
are stable equilibriums are the points $O_2$, $O_4$ and $O_6$ are
unstable ones. We can consider only the quasi-potentials between the
stable ones. Suppose we have $V^+(O_1,O_3)=1$, $V^+(O_3,O_1)=2$,
$V^+(O_1,O_5)=6$, $V^+(O_5,O_1)=7$, $V^+(O_5,O_3)=3$,
$V^+(O_3,O_5)=4$.

We are concerned with the limit $\lim\li_{\ve\da
0}u^\ve(y_1,y_2,T(\ve))$ for $T(\ve)\asymp \exp(\frac{\lb}{\ve^2})$.
Starting from the initial point $(y_1,y_2)$, we suppose that we are
attracted to $O_1$ first. By calculating $\min\li_{g\in
G(L\backslash \{1\})}\sum\li_{(m\ra n)\in g}V^+(O_m,O_n)=1$ we see
that over time $\exp(\frac{1}{\ve^2})$ we are going to jump to $O_3$
first. We then calculate $\min\li_{g\in G(L\backslash
\{3\})}\sum\li_{(m\ra n)\in g}V^+(O_m,O_n)=2$ and we see that over
time $\exp(\frac{2}{\ve^2})$ we will jump from $O_3$ back to $O_1$
and we form a cycle $\pi^{(1)}=\{1,3\}$ of rank $1$. We then
calculate $A(\pi^{(1)})=\min\li_{g\in G(L\backslash \pi^{(1)})}
\sum\li_{(m \ra n)\in g}V^+(O_m,O_n)=V^+(O_1,O_3)+V^+(O_3,O_5)=5$
and the first state out of cycle $\pi^{(1)}$ that we are going to
jump to is $O_5$. Within cycle $\pi^{(1)}$ we are mostly staying in
$O_3$. We calculate $C(\pi^{(1)})=5-\min\li_{i \in
\{1,3\}}\min\li_{g\in G_{\{1,3\}}\{i\}} \sum\li_{(m\ra n)\in
g}V^+(O_m,O_n)=4$. This means, that over time
$\exp(\frac{4}{\ve^2})$ we are jumping from $O_3$ to $O_5$. We then
calculate $\min\li_{g\in G(L\backslash \{5\})}\sum\li_{(m\ra n)\in
g}V^+(O_m,O_n)=3$ and we see that we are jumping from $O_5$ out to
$O_3$ in time $\exp(\frac{3}{\ve^2})$. This implies that within the
cycle $\pi^{(2)}=\{1,3,5\}$ which exhausts all $\om$-limit sets, we
are mostly staying in $\pi^{(1)}$, and within $\pi^{(1)}$ it is
$O_3$.

Our result can be summarized as

$$\lim\li_{\ve \da
0}u^\ve(y_1,y_2,T(\ve))=g(O_1) \text{ for } T(\ve)\asymp
\exp(\frac{\lb}{\ve^2}) \text{ and } 0<\lb<1 \ ;$$

$$\lim\li_{\ve \da
0}u^\ve(y_1,y_2,T(\ve))=g(O_3) \text{ for } T(\ve)\asymp
\exp(\frac{\lb}{\ve^2}) \text{ and } 1\leq\lb \ .$$

\

\textbf{Acknowledgement}: We would like to thank our advisor
Professor Mark Freidlin for posing this problem to us and for many
useful discussions.

\


\begin{thebibliography}{100}

\bibitem{[Anderson-Orey]} Anderson, R.F., Orey, S., Small random perturbations
of dynamical systems with reflecting boundary, \textit{Nagoya Math
J.}, \textbf{60} (1976), 189--216.

\bibitem{[F sublimiting]} Freidlin, M., Sublimiting Distributions and
Stabilization of Solutions of Parabolic Equations with a Small
Parameter, \textit{Soviet Math. Dokl.}, \textbf{235}, 5, 1042--1045,
1977.

\bibitem{[F red book]} Freidlin, M., \textit{Functional integration and partial
differential equations}, Annals of Mathematical Studies, Princeton
University Press, 1985.

\bibitem{[FW 1969]} Freidlin, M., Wentzell, A., On small random
perturbations of dynamical systems, \textit{Russ. Math. Surv.},
\textbf{25} (1970), No.1, 1--56.

\bibitem{[FW book]} Freidlin, M., Wentzell, A., \textit{Random perturbations
of dynamical systems}, Second Edition, Springer, 1998.

\bibitem{[FZ]} Freidlin, M., Zhivoglyadova, L.,
Boundary value problems with a small parameter for a diffusion
process with reflection, \textit{Russ. Math. Surv.}, \textbf{31}
(1976), No.5, 241--242 (in Russian).

\bibitem{[Gikman-Skrokhod]} Gikhman, I., Skorokhod, A., \textit{The Theory of Stochastic Processes},
III, Springer, 1979.

\bibitem{[Ito-McKean]} It\^{o}, K., McKean, H.P. Jr., \textit{Diffusion
processes and their sample paths}, Springer, 1974.

\bibitem{[Sato-Tanaka]} Sato, K., Tanaka, H., Local times on the
boundary for multidimensional reflecting diffusion, \textit{Proc.
Japan Acad.}, \textbf{38}, 10 (1962), 699--702.

\bibitem{[Sato-Ueno]} Sato, K., Ueno, T., Multidimensional diffusion
and Markov processes on the boundary, \textit{J. Math. Kyoto U.},
\textbf{4} (1965), 529--605.

\bibitem{[Watanabe]} Watanabe, S., On stochastic differential
equations for multidimensional diffucion processes with boundary, I,
II, \textit{J. Math. Kyoto U.}, \textbf{11} (1971), 169--180,
545--551.

\end{thebibliography}
\end{document}